\documentclass[letterpaper,11pt]{article}

\usepackage{graphicx,color,epsfig}
\usepackage{amsfonts,amsmath,amssymb,amsthm}
\usepackage{url}
\usepackage{comment}
\usepackage{wrapfig}
\usepackage{hyperref}
\usepackage[margin=1in]{geometry}
\usepackage{authblk}

\usepackage{hyperref}
\usepackage{lmodern}
\usepackage{pslatex}
\usepackage{amsmath, amsthm, amssymb}
\usepackage{mathrsfs}
\usepackage{commath}
\usepackage{program}
\usepackage{algorithm}
\usepackage{algorithmicx}
\usepackage{algpseudocode}
\usepackage{indentfirst}
\usepackage{flexisym}
\newtheorem{theorem}{Theorem}
\newtheorem{proposition}{Proposition}
\newtheorem{lemma}{Lemma}

\newtheorem{remark}{Remark}
\title{Functional Weak Limit of Random Walks in Cooling Random Environment}

\author[]{Yongjia Xie}
\affil[]{Purdue University}

\date{\today}

\providecommand{\keywords}[1]
{
  \small	
  \textbf{\textit{Keywords---}} #1
}
\begin{document}
\maketitle

\begin{abstract}
  We prove an annealed weak limit of the trajectory of the random walks in cooling random environment (RWCRE) under both slow (polynomial) and fast (exponential) cooling. We identify the weak limit when the underlying static environment is recurrent (Sinai's model). Avena and den Hollander have previously proved
  a Gaussian limiting distribution for the distribution of the endpoint of the walk.
  We find that the weak limit of the trajectory exists as a time-rescaled Brownian motion in the slow cooling case but the limit degenerates to a constant function in the fast cooling case.
\end{abstract}

\section{Introduction, Background, and Main Results}
Research on random walks in disordered environments has attracted a lot of attention by mathematicians and physicists over the last few decades.
The model of random walks in random environments (RWRE) was first studied by Solomon in \cite{Solomon}. In this model the disorder in the environment is random\footnote{A common assumption is that the randomness in the environment is i.i.d., though there also results where weaker assumptions are used instead.}, but fixed for all time by the walk.
Much of the subsequent interest in this model was driven by the fact that RWREs could exhibit a surprisingly rich array of asymptotic behaviors such as transience with asymptotically zero speed \cite{Solomon}, and limiting distributions which are non-Gaussian and have non-diffusive scaling \cite{Kozlov,Kesten,Sinai}.
These interesting phenomena can be understood as occurring because of the ``trapping'' effects of the environment. See \cite{ofer} for an overview of basic results in RWRE.

More recently, there has been interest in a generalization of RWRE called random walks in dynamic random environments (RWDRE) in which the disorder of the environment is random in both space and time. One can see that RWDRE interpolates between simple random walk (SRW) and RWRE: If the dynamics are "frozen", i.e. the environment is not changing after initial set-up, then this is simply a RWRE.
On the other hand, if the environment is space-time i.i.d.\ then it is easy to see that the distribution of the RWDRE (under the annealed measure) is the same as that of a SRW.
For RWDRE models which are between these two extremes there is an interplay between the trapping effects introduced by the randomness of the environment and the rate at which the time dynamics of the environment causes these traps to disappear.
One might expect therefore that environments with ``fast'' mixing time dynamics should have similar characteristics as a SRW (e.g. path convergence to Brownian motion) while ``slow'' mixing time dynamics might retain some of the strange behaviors of RWRE (e.g., non-Gaussian limiting distributions or transience with sublinear speed).

Many of the results thus far in RWDRE have focused on dynamic environments which are in some sense fast mixing, see \cite{MR2786643,BHT18,hollander6}.
For example, the environment may be assumed to be a Markov chain with uniformly mixing time dynamics or which satifies a Poincar\'e inequality.
A variety of approaches have been used in these papers, but in all cases one can obtain convergence to Brownian motion after centering and diffusive scaling.

Environments which are more slowly mixing present different problems as the trapping effects of the environment may possibly be stronger. Examples of environments like conservative particle systems have poor mixing rates \cite{DENHOLLANDER2014785,HDHS15}.
A particularly interesting example is the case where the dynamic environment is given by a simple symmetric exclusion process. Avena and Thomann have made conjectures based on simulations that this model can exhibit many of the same strange behaviors as that of RWRE (e.g., transience with zero speed and non-diffusive scaling). However, the results for this model have been limited to some cases where the parameters of the model are near their extremes and in these cases once again the distribution of the walk converges under diffusive scaling to a Brownian motion.
Other examples of slow mixing environments for which the RWDRE has been shown to converge to Brownian motion are \cite{MR3108811,Hilario,Huveneers}.

All the above results for RWDRE have shown limiting behavior which is like that of a SRW. Recently, however, Avena and den Hollander have introduced a new model of RWDRE, \emph{random walks in cooling random environment}(RWCRE), in which the dynamics can be slow enough that the model retains some of the strange behavior of RWRE \cite{hollander}.
In this model the environment is totally refreshed at some points called resampling times.
Results for this model have included a strong law of large numbers, a quenched large deviation principles, sufficient conditions for recurrence/transience, and  limiting distributions \cite{MR3940768, Hollander5, hollander}.
Most relevant to the results of the present paper, for certain cases of RWCRE they prove that the limiting distributions are Gaussian but with non-diffusive scalings that interpolate between the $(\log n)^2$ scaling of recurrent RWRE and the diffusive $\sqrt{n}$ scaling of SRW \cite{Hollander5, hollander}.
The main goal of this paper is to determine the appropriate limiting distributions for the path of the walk in these cases.

The paper is organized as follows. We introduce the model of one-dimensional RWCRE in Section 1.1. In Section 1.2 we review the limiting distribution results for both recurrent RWRE (Sinai's random walk) and corresponding model of RWCRE.
In Section 1.3, we give our main result, the functional weak limit under both slow (polynomial) and fast (exponential) cooling. The proof is given in Section 2.

\subsection{Random Walks in Cooling Random Environment}
We will use the same notations as in Avena and den Hollander \cite{hollander}. Let $\mathbb{N}_0=\mathbb{N}\cup \{0\}$. The classical one-dimensional random walk in random environment (RWRE) is defined as follows. Let $\omega =\{\omega(x):x\in\mathbb{Z}\}$ be an i.i.d. sequence with probability distribution
\begin{equation}
\mu=\alpha^{\mathbb{Z}}
\end{equation}
for some probability distribution $\alpha$ on (0,1).  The random walk in the \emph{spatial} environment $\omega$ is the Markov process $Z=(Z_n)_{n\in\mathbb{N}_0}$ starting at $Z_0=0$ with transition probabilities
\begin{equation}
P^{\omega}(Z_{n+1}=x+e|Z_n=x)=
\begin{cases}
\omega(x),~~~~~~~~ \text{if}~ e=1,\\1-\omega(x),~~ \text{if}~ e=-1,
\end{cases}
~~~~~~~n\in\mathbb{N}_0.
\end{equation}
The properties of $Z$ are well understood, both under the \emph{quenched} law $P^\omega(\cdot)$ and the \emph{annealed} law
\begin{equation}
\mathbb{P}_\mu(\cdot)=\int_{(0,1)^{\mathbb{Z}}}P^\omega(\cdot)\mu(d\omega).
\end{equation}
\bigskip

The random walk in cooling random environment (RWCRE) is a model where $\omega$ is updated along a growing sequence of determined times. Let $\tau:~~\mathbb{N}_0~\rightarrow~\mathbb{N}_0$ be a strictly increasing map such that $\tau(0)=0$ and $\tau(k)\geq k$ for $k\in \mathbb{N}$. Define a sequence of random environments $\Omega=(\omega_n)_{n\in\mathbb{N}_0}$ as follows: At each time $\tau(k)$, $k\in\mathbb{N}_0$, the environment $\omega_{\tau(k)}$ is freshly resampled from $\mu=\alpha^{\mathbb{Z}}$ and does not change during the time interval $[\tau(k),\tau(k+1))$. That is, $\omega_n = \omega_{\tau(k)}$ where $k$ is such that $\tau(k) \leq n < \tau(k+1)$. The random walk in the \emph{space-time} environment $\Omega$ is the Markov process $X=(X_n)_{n\in\mathbb{N}_0}$ starting at $X_0=0$ with transition probabilities
\begin{equation}
P^{\Omega,\tau}(X_{n+1}=x+e|X_n=x)=
\begin{cases}
\omega_n(x),~~~~~~~~ \text{if}~ e=1,\\1-\omega_n(x),~~ \text{if}~ e=-1,
\end{cases}
~~~~~~~n\in\mathbb{N}_0.
\end{equation}
We call $X$ the \emph{random walk in cooling random environment} with \emph{resampling rule} $\alpha$ and \emph{cooling rule} $\tau$.
The distribution $P^{\Omega,\tau}$ of the random walk for a given space time environment is called the \emph{quenched law}.
The \emph{annealed law} of the walk $\{X_n\}_{n\geq 0}$ is obtained by averaging the quenched with respect to the distribution $\mathbb{Q} = \mathbb{Q}_{\alpha,\tau}$ on $\Omega$.
\begin{equation}
\mathbb{P}^{\tau}(\cdot)=\int_{{((0,1)^{\mathbb{Z}})}^{\mathbb{N}_0}}P^{\Omega,\tau}(\cdot)\mathbb{Q}(d\omega),
\end{equation}
\subsection{Slow and Fast Cooling: Gaussian Fluctuations for Recurrent RWRE}
In Solomon's seminal paper \cite{Solomon}, he showed that the recurrence/transience of a RWRE is determined by the sign of $E_\alpha[\log\rho(0)]$, where
\begin{equation}
\rho(0)=\frac{1-\omega(0)}{\omega(0)}
\end{equation}
and $E_\alpha[\cdot]$ denotes expectations with respect to the measure $\alpha$. In particular, if $E_\alpha[\log \rho(0)] = 0$ then the RWRE is recurrent.
Subsequently, the scaling limit in the recurrent case was identified by Sinai \cite{Sinai} and the explicit form of the limiting distribution by Kesten \cite{Kesten}. Moreover, it was shown by Avena and den Hollander \cite{hollander} that the convergence also holds in $L^p$.The next proposition summarises their results.
\begin{proposition}\label{prop:Sinai}
\emph{{\textbf{[Scaling limit RWRE: recurrent case]}}} Let $\alpha$ be any probability distribution on (0,1) satisfying $E(\log\rho(0))=0$ and $\sigma_\mu^2=E[\log^2\rho(0)]\in(0,\infty)$. Then, under the annealed law $\mathbb{P}_\mu$, the sequence of random variables
\begin{equation}
\frac{Z_n}{\sigma_\mu^2\log^2n},~~~n\in\mathbb{N},
\end{equation}
converges in distribution and in $L^p$ to a random variable $V$ on $\mathbb{R}$ that is independent of $\alpha$. The law of $V$ has a density $p(x)$, $x\in\mathbb{R}$, with respect to the Lebesgue measure that is given by
\begin{equation}\label{Vdensity}
p(x)=\frac{2}{\pi}\sum_{k\in\mathbb{N}_0}\frac{(-1)^k}{2k+1}\exp\left[-\frac{(2k+1)^2\pi^2}{8}|x|\right],~~~~x\in\mathbb{R}.
\end{equation}
\end{proposition}
In their initial paper on RWCRE Avena and den Hollander introduced several kinds of cooling regimes that are interesting to research. For RWCRE in our paper, following their works, we focus on two kinds of growth regimes for $\tau(k)$. Let $T_k=\tau(k)-\tau(k-1)$,\\
\noindent (R1) \emph{Slow cooling:} $T_k\sim \beta Bk^{\beta-1}$, for some $B\in(0,\infty)$ and $\beta\in(1,\infty)$.\\
\noindent (R2) \emph{Fast cooling:} $\log T_k\sim Ck$, for some $C\in(0,\infty)$.

When the distribution $\alpha$ is as in Proposition \ref{prop:Sinai}, Avena and den Hollander \cite{hollander} proved a limiting distribution for the walk under both the fast and slow cooling regimes. Later in \cite{Hollander5} they strengthened this to $L^p$ convergence.  The following proposition summaries their results. Note that here and throughout the remainder of the paper we will use $\mathcal{N}(\mu,\sigma^2)$ to denote a Gaussian random variable with mean $\mu$ and variance $\sigma^2$.
\begin{proposition}\emph{\textbf{[Slow and fast cooling: Gaussian fluctuations for recurrent RWRE]}}
Let $\alpha$ be as in Proposition \ref{prop:Sinai}. In regime (R1) and (R2), under the annealed law $\mathbb{P}$,
\begin{equation}
\frac{X_n-\mathbb{E}(X_n)}{\sqrt{\chi_n(\tau)}}\rightarrow^{L^p} \mathcal{N}(0,1),
\end{equation}
where
\begin{equation}
\chi_n(\tau)=
\begin{cases}
(\sigma_\mu^2\sigma_V)^2(\frac{\beta-1}{\beta})^4(\frac{n}{B})^{\frac{1}{\beta}}\log^4 n,~~~in~regime~\emph{(R1)},\\
(\sigma_\mu^2\sigma_V)^2(\frac{1}{5C^5})\log^5 n,~~~~~~~~~~~~~~in~regime~\emph{(R2)},
\end{cases}
\end{equation}
with $\sigma_\mu^2$ the variance of the random variable $\log\rho(0)$ and $\sigma_V^2$ the variance of the random variable with density \eqref{Vdensity}.
Moreover, in (R2) the centering part can be removed. That is,
\begin{equation}
    \frac{X_n}{\sqrt{\chi_n(\tau)}}\rightarrow^{L^p}\mathcal{N}(0,1).
\end{equation}

\end{proposition}

\begin{remark}
In the most recent work \cite{Hollander5}, the authors have studied more general cooling regimes and have found their limiting behavior. In fact, despite the sequence being always tight, depending on the relative variance weight, the centered walk may not always converge. In short, relative variance weight describes how significant the variance of the walk in a single cooling interval over the variance of $X_n$. The results (Theorem 1.9 and Corollary 1.10 in \cite{Hollander5}) showed that
for general cooling sequences there might be
no limiting distribution for $X_{n}/\sqrt{\mathbb{V}ar(X_n)}$, but that one can identify a class of limiting distributions along subsequences which are mixtures of Kesten's distribution and standard Gaussian.
See Examples 5 and 6 in \cite{Hollander5} for more details.
\end{remark}

\subsection{Functional Weak Limit under the Slow and Fast Cooling}
In this section we will introduce our main results of the weak limit of $(\tilde X_k/\sqrt{\chi_n(\tau)},~k=1,2,...,n)$ where $\tilde X_k=X_k(\omega)-\mathbb{E}(X_k)$ is the centered walk\footnotemark of $X_k$ under both polynomial and exponential cooling. \footnotetext{All the \~{} signs in our paper mean the centered random variable under the annealed measure.}Since the walk $(\tilde X_k,~k=1,2,...,n)$ is a discrete time random walk and we are considering the scaled (under both time and space parameters) weak limit of it, it is reasonable to make this discrete-time-walk a continuous random walk $X_t^n$ within the time interval $t\in[0,1]$. The simplest way to solve this is to make the process piecewise linear. To this end, define
\begin{equation}\label{1}
X_t^n(\omega)=\frac{1}{\sqrt{\chi_n(\tau)}}\tilde X_{\lfloor tn\rfloor}(\omega)+(tn-\lfloor tn\rfloor)\frac{1}{\sqrt{\chi_n(\tau)}}(\tilde X_{\lfloor tn+1\rfloor}(\omega)-\tilde X_{\lfloor tn\rfloor}(\omega)).
\end{equation}
Obviously $X_t^n$ is a random function in $C[0,1]$, the space of continuous functions on $[0,1]$, equipped with the uniform topology. 
The main results are stated as follows.
\begin{theorem}\emph{\textbf{[Slow cooling: Functional weak limit for recurrent RWRE]}}\label{slow}
Let $\alpha$ be as in Proposition 1. In regime (R1), $X_t^n$ given in \eqref{1}. Under the annealed law $\mathbb{P}$,
\begin{equation}\label{functional limit polynomial}
(X_t^n,~t\in[0,1])~\Rightarrow_n~(B_{t^{1/\beta}},~t\in[0,1]),~~~~in~regime (R1),
\end{equation}
where $(B_t,t\in[0,1])$ is a standard Brownian motion. The limit in the right hand side means a time-scaled Brownian motion. The convergence in law holds in the uniform topology on $C[0,1]$.
\end{theorem}

In the exponential cooling case, the result is slightly different. The functional weak limit of $X_t^n$ is a random constant function and the law of the random constant is a standard Gaussian distribution.

\begin{theorem}\emph{\textbf{[Fast cooling: Functional weak limit for recurrent RWRE]}}\label{fast}
Let $\alpha$ be as in Proposition 1. In regime (R2), $X_t^n$ given in \eqref{1}. Under the annealed law $\mathbb{P}$, for any $a\in(0,1]$,
\begin{equation}\label{functional limit exponential}
(X_t^n,~t\in[a,1])~\Rightarrow_n~(N_t,~t\in[a,1]),~~~~~~~in~regime (R2),
\end{equation}
where $N_t=N$ for all $t\in[a,1]$ and $N\sim \mathcal{N}(0,1)$. The convergence in law holds in the uniform topology on $C[a,1]$.
\end{theorem}

\begin{remark}
In Theorem \ref{fast} the convergence holds in space $C[a,1]$ for any $a\in(0,1]$. In fact, if we want to extend the convergence to then entire time interval $[0,1]$ then neither continuous function space $C[0,1]$ nor the C\`adl\`ag function space $D[0,1]$ (with the Skorohod topology) will be sufficient since the sequence is not tight in either space.
Moreover, one can guess the limiting process on $[0,1]$ should be 0 when $t=0$ and $N_t$ for $t\in(0,1]$,
which is not a C\`adl\`ag function. So if we want to extend the convergence to a function space on $[0,1]$ then a wider space would be required, e.g. $L^p[0,1]$, together with a corresponding topology where the weak convergence holds.
\end{remark}

\section{Proof of the Theorem}
We begin by noting the following useful decomposition property of RWCRE. Let
\begin{equation}
k(n)=\max\{k\in\mathbb{N}:\tau(k)\leq n\}
\end{equation}
be the number of resamplings of the environment prior to time $n$. It's easy to see $k(n)\sim(n/B)^{1/\beta}$ in (R1) and $k(n)\sim (1/C)\log n$ in (R2). Furthermore, $X_n$ has a decomposition that will be very useful in the following proof of the theorems,
\begin{equation}\label{9}
X_n=\sum_{j=1}^{k(n)}Y_j+\bar Y_{n},
\end{equation}
where $Y_j=X_{\tau(j)}-X_{\tau(j-1)}, j=1,2,..,k(n)$, $\bar Y_{n}=X_n-X_{\tau(k(n))}$. A simple fact is that all terms in (\ref{9}) are independent under the annealed measure. Moreover, under the annealed measure, $Y_j$ has the same distribution as $Z_{T_j}$ for $j\geq 1$, and $\bar Y_n$ has the same distribution as $Z_{n-\tau(k(n))}$ for $n\geq 1$, where $\{Z_n\}_{n\geq 0}$ is a RWRE. Since we will deal with the remainder part $\bar Y_n$ throughout the proof, we will use the notation $\bar T_n=n-\tau(k(n))$ and $\bar T_n^c=\tau(k(n)+1)-n$.
\subsection{Slow Cooling}
We start by finding the weak limit of the finite dimensional random vector $(X_{t_1}^n,X_{t_2}^n,...,X_{t_k}^n)$. To start with, we will prove the weak convergence under the case $k=2$, i.e. the weak limit of $(X_t^n,X_s^n)$ for $0\leq t< s\leq 1$. By \cite{hollander}, $\tilde X_{\lfloor tn\rfloor}/\sqrt{\chi_{\lfloor tn\rfloor}(\tau)}\Rightarrow_n \mathcal{N}(0,1)$. Obviously $\lim \frac{\chi_{\lfloor tn\rfloor}(\tau)}{\chi_n(\tau)}=t^{1/\beta}$, so $\tilde X_{\lfloor tn\rfloor}/\sqrt{\chi_n(\tau)}\Rightarrow_n \mathcal{N}(0,t^{1/\beta})$. If $\psi_{n,t}$ is the rightmost term in \eqref{1}, then $\psi_{n,t}\Rightarrow_n 0$ by the fact that all the numerators are bounded but $\chi_n(\tau)$ goes to infinity. We have
\begin{equation}\label{2}
(X_t^n,X_s^n-X_t^n)=\frac{1}{\sqrt{\chi_n(\tau)}}(\tilde X_{\lfloor tn\rfloor}(\omega),~\tilde X_{\lfloor sn\rfloor}(\omega)-\tilde X_{\lfloor tn\rfloor}(\omega))+(\psi_{n,t},~\psi_{n,s}-\psi_{n,t}).
\end{equation}
To find the weak limit of $(\tilde X_{\lfloor sn\rfloor}-\tilde X_{\lfloor tn\rfloor})/\sqrt{\chi_n(\tau)}$,
we will follow the approach of \cite{hollander} in using the following Lyapunov condition.
\begin{lemma}\emph{\textbf{(Lyapunov condition, Petrov \cite{Petrov})}}\\
Let $U=(U_k)_{k\in \mathbb{N}}$ be a sequence of independent random variables (at least one of which has a non-degenerate distribution). Let $m_k=E(U_k)$ and $\sigma_k^2=Var(U_k)$. Define
\begin{equation}
\chi_n=\sum_{k=1}^n \sigma_k^2.
\end{equation}
Then the Lyapunov condition
\begin{equation}\label{lyapunov}
\lim_{n\rightarrow\infty}\frac{1}{\chi_n^{p/2}}\sum_{k=1}^nE(|U_k-m_k|^p)=0,
\end{equation}
for some $p>2$
implies that
\begin{equation}
\frac{1}{\chi_n}\sum_{k=1}^n(U_k-m_k)\Rightarrow_n~\mathcal{N}(0,1).
\end{equation}.

\end{lemma}
Recall $X_n$ has the decomposition
\begin{equation}
X_n=\sum_{j=1}^{k(n)}Y_j+\bar Y_{n}.
\end{equation}
Define the variance of $X_{\lfloor sn\rfloor}-X_{\lfloor tn\rfloor}$ (which is also the variance of $\tilde X_{\lfloor sn\rfloor}-\tilde X_{\lfloor tn\rfloor}$) for any $s<t$ and $n$ large enough
\begin{equation}\label{Varst}
\chi_n^{t,s}(\tau)=\sum_{j=k(\lfloor tn\rfloor)+2}^{k(\lfloor sn\rfloor)}\mathbb{V}ar(Y_j)+\mathbb{V}ar(\bar Y_{\lfloor sn\rfloor})+\mathbb{V}ar(\bar Y_{\lfloor tn\rfloor}^c)
\end{equation}
where $\bar Y_n^c=X_{\tau(k(n)+1)}-X_n$.
Recall that $\tilde Y_j=Y_j-\mathbb{E}(Y_j)$, $\tilde{\bar{Y}}_n=\bar Y_n-\mathbb{E}(\bar Y_n)$, and $\tilde{\bar{Y}}_n^c=\bar Y_n^c-\mathbb{E}(\bar Y_n^c)$. For $p>2$, let
\begin{equation}\label{pmst}
\chi_n^{t,s}(\tau;p)=\sum_{j=k(\lfloor tn\rfloor)+2}^{k(\lfloor sn\rfloor)}\mathbb{E}\left(|\tilde Y_j|^p\right)+\mathbb{E}\left(|\tilde {\bar Y}_{\lfloor sn\rfloor}|^p\right)+\mathbb{E}\left(|\tilde{\bar Y}_{\lfloor tn\rfloor}^c|^p\right).
\end{equation}
Since $Y_j$ has the same distribution as $Z_{T_j}$, then by Proposition 4 in \cite{hollander} the following two asymptotic estimates hold as $j\rightarrow\infty$.
\begin{equation}\label{estimation}
\mathbb{V}ar(Y_j)\sim (\sigma_\mu^2\sigma_V)^2\log^4T_j,~~~~\mathbb{E}\left(|\tilde Y_j|^p\right)=O(\log^{2p}T_j),~~p>2.
\end{equation}
Applying these to \eqref{Varst} and \eqref{pmst} we obtain
\begin{equation}\label{23}
\begin{split}
&\sum_{j=k(\lfloor tn\rfloor)+2}^{k(\lfloor sn\rfloor)}\mathbb{V}ar(Y_j)\sim (\sigma_\mu^2\sigma_V)^2\sum_{j=k(\lfloor tn\rfloor)+2}^{k(\lfloor sn\rfloor)}\log^4T_j,\\
&\sum_{j=k(\lfloor tn\rfloor)+2}^{k(\lfloor sn\rfloor)}\mathbb{E}\left(|\tilde Y_j|^p\right)=O\left(\sum_{j=k(\lfloor tn\rfloor)+2}^{k(\lfloor sn\rfloor)}\log^{2p}T_j\right).
\end{split}
\end{equation}
Moreover, using that $\sum_{j=1}^k\log^{2p}j\sim \int_1^k \log^{2p}xdx\sim k\log^{2p}k$ for all $p\geq 2$ and that $k(n) \sim (n/B)^{1/\beta}$, we have
\begin{equation}\label{30}
\begin{split}
(\sigma_\mu^2\sigma_V)^2\sum_{j=k(\lfloor tn\rfloor)+2}^{k(\lfloor sn\rfloor)}\log^4T_j&\sim (\sigma_\mu^2\sigma_V)^2(\beta-1)^4\left[(\frac{sn}{B})^{1/\beta}\log^4\left((\frac{sn}{B})^{1/\beta}\right)-(\frac{tn}{B})^{1/\beta}\log^4\left((\frac{tn}{B})^{1/\beta}\right)\right]\\ &\sim (\sigma_\mu^2\sigma_V)^2(\beta-1)^4(\frac{n}{B})^{1/\beta}(\frac{1}{\beta})^4(s^{1/\beta}-t^{1/\beta})\log^4n\\
&=\chi_n(\tau)\left(s^{\frac{1}{\beta}}-t^{\frac{1}{\beta}}\right),\\
\text{and} \qquad \sum_{j=k(\lfloor tn\rfloor)+2}^{k(\lfloor sn\rfloor)}\log^{2p}T_j&\sim (\sigma_\mu^2\sigma_V)^2(\beta-1)^{2p}(\frac{n}{B})^{1/\beta}(\frac{1}{\beta})^{2p}(s^{1/\beta}-t^{1/\beta})\log^{2p}n\\
&=\chi_n^{\frac{p}{2}}(\tau)\left(s^{\frac{1}{\beta}}-t^{\frac{1}{\beta}}\right)\left(\frac{n}{B}\right)^{\frac{2-p}{\beta p}},~~~~p>2.
\end{split}
\end{equation}
Since $\bar Y_n$ has the same distribution as $Z_{\bar T_n}$, we can again use Proposition 4 in \cite{hollander} to obtain that there exists $C^{(2)}>0,~~C^{(p)}>0$, such that
\begin{equation}\label{26}
\mathbb{V}ar\left(\bar Y_{{n}}\right)\leq C^{(2)}\log^4\bar T_{n},~~~~
\mathbb{E}\left(|\tilde{\bar Y}_{n}|^p\right)\leq C^{(p)}\log^{2p}\bar T_{n}.
\end{equation}
These upper bounds will be used to control $\mathbb{V}ar(\bar Y_{n}^c)$ and $\mathbb{E}\left(|\bar Y_{n}^c-\mathbb{E}(\bar Y_{n}^c)|^p\right)$. For $n$ large enough,
\begin{equation}\label{25}
\begin{split}
&\mathbb{V}ar(\bar Y_{n}^c)=\mathbb{V}ar(Y_{k(n)+1}-\bar Y_{n})\leq2\mathbb{V}ar(Y_{k(n)+1})+2\mathbb{V}ar(\bar Y_n)\leq 4\left[(\sigma_\mu^2\sigma_V)^2+C^{(2)}\right]\log^4 T_{k(n)+1},\\
&\mathbb{E}\left(|\tilde{\bar Y}_{n}^c|^p\right)=\mathbb{E}\left(|\tilde Y_{k(n)+1}-\tilde{\bar Y}_n|^p\right)\leq 2^{p-1}\left[\mathbb{E}\left(|\tilde Y_{k(n)+1}|^p\right)+\mathbb{E}\left(|\tilde{\bar Y}_n|^p\right)\right]=O(\log^{2p}T_{k(n)+1}).
\end{split}
\end{equation}
From (\ref{26}) and (\ref{25}),
\begin{equation}\label{29}
\begin{split}
&\mathbb{V}ar(\bar Y_{\lfloor sn\rfloor})+\mathbb{V}ar(\bar Y_{\lfloor tn\rfloor}^c)\leq C^{(2)}\log^4\bar T_{\lfloor sn\rfloor}+4\left[(\sigma_\mu^2\sigma_V)^2+C^{(2)}\right]\log^4 T_{k(\lfloor tn\rfloor)+1}=O(\log^4 n),\\
&\mathbb{E}\left(|\tilde {\bar Y}_{\lfloor sn\rfloor}|^p\right)+\mathbb{E}\left(|\tilde{\bar Y}_{\lfloor tn\rfloor}^c|^p\right)\leq
C^{(p)}\log^{2p}\bar T_{\lfloor sn\rfloor}+O(\log^{2p}T_{k(\lfloor tn\rfloor)+1})=O(\log^{2p}n).
\end{split}
\end{equation}
By \eqref{30} and \eqref{29}, we can therefore give the asymptotic of $\chi_n^{t,s}(\tau)$ and $\chi_n^{t,s}(\tau;p)$,
\begin{equation}
\begin{split}
\chi_n^{t,s}(\tau)&\sim \chi_n(\tau)\left(s^{\frac{1}{\beta}}-t^{\frac{1}{\beta}}\right),\\
\chi_n^{t,s}(\tau;p)&=O\left(\chi_n^{\frac{p}{2}}(\tau)\left(s^{\frac{1}{\beta}}-t^{\frac{1}{\beta}}\right)\left(\frac{n}{B}\right)^{\frac{2-p}{\beta p}}\right),~~~~p>2.
\end{split}
\end{equation}
From these asymptotics it is easy to check that the Lyapunov condition holds, and thus
\begin{equation}
\frac{\tilde X_{\lfloor sn\rfloor}-\tilde X_{\lfloor tn\rfloor}}{\sqrt{\chi_n(\tau)}}\Rightarrow_n \mathcal{N}(0,s^{1/\beta}-t^{1/\beta}).
\end{equation}
In order to prove the vector $(X_t^n,X_s^n-X_t^n)$ converges to a 2-d Gaussian vector with independent components, it suffices to show that any linear combination of $X_t^n$ and $X_s^n-X_t^n$ converges to the corresponding linear combination of the components of the limiting Gaussian vector. To this end, the proof is quite similar to what we did above: Decompose $\lambda X_{\lfloor tn\rfloor}+\mu (X_{\lfloor sn\rfloor}-X_{\lfloor tn\rfloor})$ into independent sums and check the Lyapunov condition (\ref{lyapunov}). Notice that
\begin{equation}
\lambda  X_{\lfloor tn\rfloor}+\mu ( X_{\lfloor sn\rfloor}- X_{\lfloor tn\rfloor})=\lambda\sum_{j=1}^{k(\lfloor tn\rfloor)}Y_j+\left(\lambda\bar Y_{\lfloor tn\rfloor}+\mu\bar Y_{\lfloor tn\rfloor}^c\right)+\mu\sum_{j=k(\lfloor tn\rfloor)+2}^{k(\lfloor sn\rfloor)}Y_j+\mu\bar Y_{\lfloor sn\rfloor}.
\end{equation}
The key point to the proof is the expressions of the variance  of $\lambda  X_{\lfloor tn\rfloor}+\mu ( X_{\lfloor sn\rfloor}- X_{\lfloor tn\rfloor})$
\begin{equation}\label{3}
\begin{split}
\mathbb{V}ar\left(\lambda X_{\lfloor tn\rfloor}+\mu ( X_{\lfloor sn\rfloor}- X_{\lfloor tn\rfloor})\right)=\lambda^2\sum_{j=1}^{k(\lfloor tn\rfloor)}\mathbb{V}ar(Y_j)+\mu^2\sum_{j=k(\lfloor tn\rfloor)+2}^{k(\lfloor sn\rfloor)}\mathbb{V}ar(Y_j)\\ +\mu^2\mathbb{V}ar(\bar Y_{\lfloor sn\rfloor})+
\mathbb{V}ar\left(\lambda\bar Y_{\lfloor tn\rfloor}+\mu\bar Y_{\lfloor tn\rfloor}^c\right),
\end{split}
\end{equation}
and the sum of centered $p$-th moments of the independent components in the above decomposition,
\begin{equation}\label{4}
\begin{split}
\lambda^p\sum_{j=1}^{k(\lfloor tn\rfloor)}\mathbb{E}\left(|\tilde Y_j|^p\right)+\mu^p\sum_{j=k(\lfloor tn\rfloor)+2}^{k(\lfloor sn\rfloor)}\mathbb{E}\left(|\tilde Y_j|^p\right)+\mu^p\mathbb{E}\left(|\tilde{\bar Y}_{\lfloor sn\rfloor}|^p\right)+
\mathbb{E}\left(|\lambda\tilde{\bar Y}_{\lfloor tn\rfloor}+\mu\tilde{\bar Y}_{\lfloor tn\rfloor}^c|^p\right).
\end{split}
\end{equation}
The last term in each expression above cannot be separated into two parts because those two random variables are not independent under the annealed measure. But still, we can estimate the last term by the fact that $Var(X+Y)\leq 2(Var(X)+Var(Y))$ (and similarly, $E(|X+Y|^p)\leq 2^{p-1}(E|X|^p+E|Y|^p)$ for the $p$-th moment) for any two random variables $X$ and $Y$. Thus, with the same approach, the last two terms in \eqref{3} and \eqref{4} will be dominated by the first two sums.
Moreover, the asymptotics of the first two sums in \eqref{3} and \eqref{4} can be obtained using the same methods as in the first part of the proof above.

The result is for any $\lambda>0$, $\mu>0$, $\lambda X_t^n+\mu (X_s^n-X_t^n)$ converges weakly to $\mathcal{N}(0, \lambda^2 t^{1/\beta}+\mu^2 (s^{1/\beta}-t^{1/\beta}))$. This also reveals the independence of the coordinates of the limit random vector, i.e.
\begin{equation}
(X_t^n,X_s^n-X_t^n)\Rightarrow_n (N_1,N_2),
\end{equation}
where $(N_1,N_2)$ is a Gaussian vector with mean $(0,0)$ and variance $(t^{1/\beta},s^{1/\beta}-t^{1/\beta})$, also $N_1$ and $N_2$ are independent.

It is natural to extend the weak convergence of 2-dimension vector into finite dimension vector $(X_{t_1}^n,X_{t_2}^n,...,X_{t_k}^n)$, $0\leq t_1<t_2<...<t_k\leq1$ i.e.
\begin{equation}
(X_{t_1}^n,X_{t_2}^n,...,X_{t_k}^n)\Rightarrow_n (B_{t_1^{1/\beta}},B_{t_2^{1/\beta}},...,B_{t_k^{1/\beta}}),
\end{equation}
where $(B_t,t\in[0,1])$ is a standard Brownian motion. The proof of this statement follows the same steps as what we did in dimension 2: Decompose $\sum_{i=1}^k\lambda_i(X_{\lfloor t_in\rfloor}-X_{\lfloor t_{i-1}n\rfloor})$ into independent sums  where $t_0=0$. Then take the variance and the the sum of centered $p-$th moment of the independent components of the decomposition to check the Lyapunov condition \eqref{lyapunov}. The decomposition is
\begin{equation}
\begin{split}
    \sum_{i=1}^k\lambda_i(X_{\lfloor t_in\rfloor}-X_{\lfloor t_{i-1}n\rfloor})&=
    \lambda_1\sum_{j=1}^{k(\lfloor t_1n\rfloor)}Y_j
+\sum_{i=2}^{k}\lambda_i\sum_{j=k(\lfloor t_{i-1}n+2\rfloor)}^{k(\lfloor t_in\rfloor)}Y_j\\
&+\left[\sum_{i=1}^{k-1}(\lambda_i\bar Y_{\lfloor t_in\rfloor}+\lambda_{i+1}\bar Y_{\lfloor t_in\rfloor}^c)+\lambda_k\bar Y_{\lfloor t_kn\rfloor}\right].
\end{split}
\end{equation}
So the variance and the sum of centered $p-$th moment of the independent components above are
\begin{equation}
\begin{split}
\lambda_1^2\sum_{j=1}^{k(\lfloor t_1n\rfloor)}\mathbb{V}ar(Y_j)
&+\sum_{i=2}^{k}\lambda_i^2\sum_{j=k(\lfloor t_{i-1}n+2\rfloor)}^{k(\lfloor t_in\rfloor)}\mathbb{V}ar(Y_j)
\\ &+\left[\sum_{i=1}^{k-1}\mathbb{V}ar\left(\lambda_i\bar Y_{\lfloor t_in\rfloor}+\lambda_{i+1}\bar Y_{\lfloor t_in\rfloor}^c\right)+\lambda_k^2\mathbb{V}ar(\bar Y_{\lfloor t_kn\rfloor})\right]
\end{split}
\end{equation}
and
\begin{equation}
\begin{split}
\lambda_1^p\sum_{j=1}^{k(\lfloor t_1n\rfloor)}\mathbb{E}(|\tilde{Y}_j|^p)
&+\sum_{i=2}^{k}\lambda_i^p\sum_{j=k(\lfloor t_{i-1}n+2\rfloor)}^{k(\lfloor t_in\rfloor)}\mathbb{E}(|\tilde{Y}_j|^p)
\\ &+\left[\sum_{i=1}^{k-1}\mathbb{E}\left(|\lambda_i\tilde{\bar Y}_{\lfloor t_in\rfloor}+\lambda_{i+1}\tilde{\bar Y}_{\lfloor t_in\rfloor}^c|^p\right)+\lambda_k^p\mathbb{E}(|\tilde{\bar Y}_{\lfloor t_kn\rfloor}|^p)\right].
\end{split}
\end{equation}
All the terms in the big brackets are dominated by sums to the left of the brackets. To check the Lyapunov condition holds in this case is nothing new but repeat our works (\ref{23}) and (\ref{30}). The details are tedious and we omit them in our paper.
\bigskip

To complete the proof of the theorem under the slow cooling case, the tightness of the sequence $X^n$ is needed. To this end, by Theorems 7.3 and 7.4 in \cite{billingsley} it is enough to show that
for any $\epsilon >0,~~\eta>0$, $\exists \delta>0$ and a sequence of numbers \{$t_i$\}, where $0=t_0<t_1<...<t_v=1$, s.t.
\begin{equation}
\min_{1<i<v}(t_i-t_{i-1})\geq \delta,
\end{equation}
and $\exists n_0>0$, for all $n>n_0$,
\begin{equation}\label{tightness}
\sum_{i=1}^v \mathbb{P}[\sup_{t_{i-1}\leq s \leq t_i}|X_s^n-X_{t_{i-1}}^n|\geq \epsilon]<\eta.
\end{equation}
Since $(X_t^n,~t\in [0,1])$ is the continuous process of $(\tilde X_{\lfloor tn\rfloor}/\sqrt{\chi_n(\tau)},~t\in [0,1])$, the biggest difference in the continuous process within a given interval is, up to an error smaller than $2/\sqrt{\chi_n(\tau)}$, bounded by the biggest difference in the discrete time process. Hence we can check the condition \eqref{tightness} by replacing $X_s^n, s\in[t_{i-1},t_i]$ and $X_{t_{i-1}}^n$ by $\tilde X_m/\sqrt{\chi_n(\tau)}$, $m\in [t_{i-1}n,t_in]$ and $\tilde X_{\lfloor t_{i-1}n\rfloor}/\sqrt{\chi_n(\tau)}$ separately.

Let $m$ be $|\tilde X_{m}-\tilde X_{\lfloor t_{i-1}n\rfloor}|=\sup_{s\in [t_{i-1}n,t_in]}|\tilde X_{s}-\tilde X_{\lfloor t_{i-1}n\rfloor}|$, i.e. the exact value of $s$ to make the biggest difference happens. If there are more than one candidates, choose one arbitrarily. We have the following decomposition,
\begin{equation}\label{decomposition}
\tilde X_{m}-\tilde X_{\lfloor t_{i-1}n\rfloor}=\sum_{j=\tau(k(\lfloor t_{i-1}n\rfloor)+1)}^{\tau(k(m))}\tilde Y_j+\tilde{\bar{Y}}_{m}-\tilde{\bar{Y}}_{\lfloor{t_{i-1}n}\rfloor},
\end{equation}
or just $\tilde{\bar{Y}}_{m}-\tilde{\bar{Y}}_{\lfloor{t_{i-1}n}\rfloor}$ if $k(\lfloor t_{i-1}n\rfloor)=k(m)$.

Let's deal with the decomposition above in two parts:
\begin{itemize}
  \item Given $q=\lfloor \beta \rfloor+1>1$, define the martingale $\{M_l\}$ as $M_0=0$,
      \begin{equation}M_l=\sum_{j=\tau(k(\lfloor t_{i-1}n\rfloor)+1)}^{\tau(k(\lfloor t_{i-1}n\rfloor)+l)}\tilde Y_j,~~~~l\geq1.
      \end{equation}
      Since the function $x^{2q}$ is convex, $\{M_l^{2q}\}$ is a submartingale. By Doob's Maximal Inequality \cite{Durret}, for integer $L>0$,
      \begin{equation}\label{doob}
      \mathbb{P}\left(\sup_{l\in [0,L]}\frac{|M_l|}{\sqrt{\chi_n(\tau)}}\geq \frac{\epsilon}{2}\right)\leq
       \frac{\mathbb{E}[M_L^{2q}]}{(\frac{\epsilon}{2})^{2q}\chi_n^q(\tau)}.
      \end{equation}
      To estimate the order of $\mathbb{E}[M_L^{2q}]$, notice that if we expand all the terms in $M_L^{2q}$, it is a sum that several terms in it have zero mean. So by counting the number of non-zero terms in $\mathbb{E}[M_L^{2q}]$ will give us the order of it. In fact, any term that has non-zero mean cannot have a factor $\tilde Y_j$ of order only one, i.e. either it is not divided by $\tilde Y_j$ or it is divided by $\tilde Y_j^2 $. Thus, a rough upper bound of the number of the non-zero terms in $\mathbb{E}[M_L^{2q}]$ is $\sum_{i=1}^q\binom{L}{i}i^{2q}$.
      Since $q$ is fixed, for $L$ large enough, $\sum_{i=1}^q\binom{L}{i}i^{2q}\leq q\binom{L}{q}q^{2q}$.

      For any nonzero term in the expansion of $\mathbb{E}[M_L^{2q}]$, by \eqref{estimation}, it is bounded from above by $C_0\log^{4q}n$ for some $C_0>0$ since we are dealing with the case within the interval $[0,n]$. So
      \begin{equation}\label{bound}
      \mathbb{E}[M_L^{2q}]\leq C_0q\binom{L}{q}q^{2q}\log^{4q}n\leq C_0q^{2q+1}L^q \log^{4q}n.
      \end{equation}

      Now back to the first part in \eqref{decomposition},
      \begin{equation}
      \mathbb{P}\left(\frac{|\sum_{j=\tau(k(\lfloor t_{i-1}n\rfloor)+1)}^{\tau(k(m))}\tilde Y_j|}{\sqrt{\chi_n(\tau)}}\geq \frac{\epsilon}{2}\right)\leq
      \mathbb{P}\left(\sup_{l\in [0,k(\lfloor t_in\rfloor)-k(\lfloor t_{i-1}n\rfloor)]}\frac{|M_l|}{\sqrt{\chi_n(\tau)}}\geq \frac{\epsilon}{2}\right).
      \end{equation}
      Combining with \eqref{doob} and \eqref{bound} and recalling $(k(\lfloor t_in\rfloor)-k(\lfloor t_{i-1}n\rfloor))\sim (n/B)^{1/\beta}(t_i^{1/\beta}-t_{i-1}^{1/\beta})$, we obtain that there exists $C^*>0$, depending only on $\epsilon$, such that
      \begin{equation}\label{16}
      \mathbb{P}\left(\frac{|\sum_{j=\tau(k(\lfloor t_{i-1}n\rfloor)+1)}^{\tau(k(m))}\tilde Y_j|}{\sqrt{\chi_n(\tau)}}\geq \frac{\epsilon}{2}\right)\leq C^*(t_i^{\frac{1}{\beta}}-t_{i-1}^{\frac{1}{\beta}})^q.
      \end{equation}
  \item To deal with $\tilde{\bar{Y}}_{m}$, notice that $|\tilde{\bar{Y}}_{m}|$ is bounded by the maximum of $|\bar Y_n-\mathbb{E}(\bar Y_n)|$ where $n\in [\tau(k(m)),\tau(k(m)+1)]$. Define $\tilde Y_j^*=\max_{n\in [\tau(j-1),\tau(j)]}|\bar Y_n-\mathbb{E}(\bar Y_n)|$, then $|\tilde{\bar{Y}}_{m}|\leq \tilde Y_{k(m)+1}^*$, where $k(m)$ can be from $k(\lfloor t_{i-1}n\rfloor)$ to $k(\lfloor t_{i}n\rfloor)$. Hence,
      \begin{equation}\label{17}
      \begin{split}
      \mathbb{P}\left(\frac{|\tilde{\bar{Y}}_{m}|}{\sqrt{\chi_n(\tau)}}\geq \frac{\epsilon}{4}\right)\leq &
      \mathbb{P}\left(\sup_{j\in [k(\lfloor t_{i-1}n\rfloor)+1,k(\lfloor t_{i}n\rfloor)+1]}\frac{\tilde{Y}_{j}^*}{\sqrt{\chi_n(\tau)}}\geq \frac{\epsilon}{4}\right)\\ \leq &
      \sum_{j=k(\lfloor t_{i-1}n\rfloor)+1}^{k(\lfloor t_{i}n\rfloor)+1}\mathbb{P}\left(\frac{\tilde{Y}_{j}^*}{\sqrt{\chi_n(\tau)}}\geq \frac{\epsilon}{4}\right).
      \end{split}
      \end{equation}

      Let $Y_j^*=\max_{n\in [\tau(j-1),\tau(j)]}|\bar Y_n|$, we have
      \begin{equation}\label{18}
      \tilde Y_j^*=\max_{n\in [\tau(j-1),\tau(j)]}|\bar Y_n-\mathbb{E}(\bar Y_n)|\leq \max_{n\in [\tau(j-1),\tau(j)]}|\bar Y_n|+\max_{n\in [\tau(j-1),\tau(j)]}\mathbb{E}|\bar Y_n|\leq Y_j^*+\mathbb{E}Y_j^*.
      \end{equation}
      Moreover, by the same proof of the Proposition 4 in \cite{hollander} (both $Z_n>a$ and $Z_n^*>a$ mean $T(a)<n$), for all $p>0$,
      \begin{equation}\label{19}
      \sup_{1\leq j\leq k(n)+1}\mathbb{E}\left(\frac{Y_j^*}{\log^2 n}\right)^p\leq\sup_{1\leq j\leq k(n)+1}\mathbb{E}\left(\frac{Y_j^*}{\log^2 T_j}\right)^p<\infty.
      \end{equation}

      From \eqref{18}, Chebyshev's Inequality, and \eqref{19}, there exists $C'>0$ depending only on $\epsilon$ such that
      \begin{equation}
      \mathbb{P}\left(\frac{\tilde{Y}_{j}^*}{\sqrt{\chi_n(\tau)}}\geq \frac{\epsilon}{4}\right)\leq
      \mathbb{P}\left(\frac{Y_j^*+\mathbb{E}Y_j^*}{\sqrt{\chi_n(\tau)}}\geq \frac{\epsilon}{4}\right)\leq
      \frac{\mathbb{E}\left(Y_j^*+\mathbb{E}Y_j^*\right)^4}{\left(\frac{\epsilon}{4}\right)^4\chi_n^2(\tau)}\leq C'n^{-\frac{2}{\beta}}.
      \end{equation}

      Now the upper bound of \eqref{17} is clear,
      \begin{equation}\label{21}
      \mathbb{P}\left(\frac{|\tilde{\bar{Y}}_{m}|}{\sqrt{\chi_n(\tau)}}\geq \frac{\epsilon}{4}\right)\leq \sum_{j=k(\lfloor t_{i-1}n\rfloor)+1}^{k(\lfloor t_{i}n\rfloor)+1}C'n^{-\frac{2}{\beta}}=C'\left(k(\lfloor t_in\rfloor)-k(\lfloor t_{i-1}n\rfloor)\right)n^{-\frac{2}{\beta}}.
      \end{equation}
      The right hand side goes to zero as n goes to infinity since $k(n)\sim (n/B)^{1/\beta}$.
\end{itemize}
\bigskip

Back to the tightness condition \eqref{tightness}, for any given $\epsilon>0,\eta>0$, let $\delta =1/K$, and $t_i=i/K,~i=0,1,...,K$, the positive integer $K$ to be determined. By (\ref{26}), (\ref{decomposition}), (\ref{16}), and (\ref{21}), there exists $c>0$ such that
\begin{equation}\label{finaltight}
\begin{split}
&\sum_{i=1}^K\mathbb{P}\left(\frac{\sup_{s\in [t_{i-1}n,t_in]}|\tilde X_{s}-\tilde X_{\lfloor t_{i-1}n\rfloor}|}{\sqrt{\chi_n(\tau)}}\geq \epsilon\right)=
\sum_{i=1}^K\mathbb{P}\left(\frac{|\tilde X_{m}-\tilde X_{\lfloor t_{i-1}n\rfloor}|}{\sqrt{\chi_n(\tau)}}\geq \epsilon\right)\\
\leq &\sum_{i=1}^K\left[\mathbb{P}\left(\frac{|\sum_{j=\tau(k(\lfloor t_{i-1}n\rfloor)+1)}^{\tau(k(m))}\tilde Y_j|}{\sqrt{\chi_n(\tau)}}\geq \frac{\epsilon}{2}\right)+\mathbb{P}\left(\frac{|\tilde{\bar{Y}}_{m}|}{\sqrt{\chi_n(\tau)}}\geq \frac{\epsilon}{4}\right)+\mathbb{P}\left(\frac{|\tilde{\bar{Y}}_{\lfloor t_{i-1}n\rfloor}|}{\sqrt{\chi_n(\tau)}}\geq \frac{\epsilon}{4}\right)\right]\\
\leq &\sum_{i=1}^K\left[C^*(t_i^{\frac{1}{\beta}}-t_{i-1}^{\frac{1}{\beta}})^q+C'\left(k(\lfloor t_in\rfloor)-k(\lfloor t_{i-1}n\rfloor)\right)n^{-\frac{2}{\beta}}+\frac{16\mathbb{V}ar(\bar Y_{\lfloor t_{i-1}n\rfloor})}{\epsilon^2\chi_n(\tau)}   \right]\\
\leq & C^*K\sup_{1\leq i\leq K}\left[\left(\frac{i}{K}\right)^{\frac{1}{\beta}}-\left(\frac{i-1}{K}\right)^{\frac{1}{\beta}}\right]^q+ c K n^{-\frac{1}{\beta}}\\
=&C^*K^{1-\frac{q}{\beta}}+c K n^{-\frac{1}{\beta}}.
\end{split}
\end{equation}
Since $q>\beta$, by first choosing $K$ large and then choosing $n$ large enough the above bound is less than $\eta$.
Hence the tightness condition holds, and \eqref{functional limit polynomial} is proved. \qed
\subsection{Fast cooling}
We do the proof in the same way as above: Find the order of the variance $\chi_n^{t,s}(\tau)$ then determine the limit distribution of $(\tilde X_{\lfloor sn\rfloor}-\tilde X_{\lfloor tn\rfloor})/\sqrt{\chi_n(\tau)}$. Check the tightness of the distribution of the process $(X_t^n,~t\in[a,1])$ to get the desired result.

Given $0<a\leq t\leq s\leq 1$, recall the variance of $\tilde X_{\lfloor sn\rfloor}-\tilde X_{\lfloor tn\rfloor}$ is
\begin{equation}
\chi_n^{t,s}(\tau)=\sum_{j=k(\lfloor tn\rfloor)+2}^{k(\lfloor sn\rfloor)}\mathbb{V}ar(Y_j)+\mathbb{V}ar(\bar Y_{\lfloor sn\rfloor})+\mathbb{V}ar(\bar Y_{\lfloor tn\rfloor}^c).
\end{equation}
Again by \eqref{estimation}, since $t>0$, for $n$ large enough,
\begin{equation}
\mathbb{V}ar(Y_j)\leq 2 (\sigma_\mu^2\sigma_V)^2\log^4 T_j\leq2 (\sigma_\mu^2\sigma_V)^2\log^4 n
\end{equation}
holds for $j\in[k(\lfloor tn\rfloor)+2,k(\lfloor sn\rfloor)]$. Recall the upper bound in \eqref{29}
\begin{equation}
\mathbb{V}ar(\bar Y_{\lfloor sn\rfloor})+\mathbb{V}ar(\bar Y_{\lfloor tn\rfloor}^c)\leq C^{(2)}\log^4\bar T_{\lfloor sn\rfloor}+4\left[(\sigma_\mu^2\sigma_V)^2+C^{(2)}\right]\log^4 T_{k(\lfloor tn\rfloor)+1}=O(\log^4 n).
\end{equation}
In the fast cooling case, since $k(n)\sim (1/C)\log n$, the number of terms in the sum $\sum_{j=k(\lfloor tn\rfloor)+2}^{k(\lfloor sn\rfloor)}\mathbb{V}ar(Y_j)$ is $(1/C)[\log{sn}+o(\log{sn})-\log{tn}-o(\log{tn})]=(1/C)[\log{s/t}+o(\log n)] $. Thus,
\begin{equation}
\chi_n^{t,s}(\tau)\leq 2(\sigma_\mu^2\sigma_V)^2\left[\frac{1}{C}\log\left(\frac{s}{t}\right)+o(\log n)\right]\log^4 n+O(\log^4 n).
\end{equation}
Since $\chi_n(\tau)$ is of order $\log^5 n$, it is obvious $(\tilde X_{\lfloor sn\rfloor}-\tilde X_{\lfloor tn\rfloor})/\sqrt{\chi_n(\tau)}\Rightarrow_n 0$. Moreover, notice that $\chi_{\lfloor tn\rfloor}(\tau)\sim \chi_n(\tau)$ for any $t\in[a,1]$, $\tilde X_{\lfloor tn\rfloor}/\sqrt{\chi_n(\tau)}\Rightarrow_n N$, where $N$ is the standard Gaussian random variable. Hence $(X_t^n,X_s^n-X_t^n)\Rightarrow_n (N,0)$.

Extend the weak convergence of 2-dimension vector into finite dimension vector $(X_{t_1}^n,X_{t_2}^n,...,X_{t_k}^n)$, i.e.
\begin{equation}
(X_{t_1}^n,X_{t_2}^n,...,X_{t_k}^n)\Rightarrow_n (N,N,...,N),
\end{equation}
where $t_i\in[a,1],i=1,2,...,k$, $N\sim \mathcal{N}(0,1)$.
\bigskip

To check the tightness condition, it is enough to show that (let $\delta=1$) for any $ \epsilon>0$,
\begin{equation}
\limsup_{n\rightarrow\infty}\mathbb{P}\left(\sup_{a\leq s\leq t\leq 1}|X_s^n-X_t^n|\geq \epsilon\right)=0,
\end{equation}
which is equivalent to
\begin{equation}\label{24}
\limsup_{n\rightarrow\infty}\mathbb{P}\left(\sup_{\lfloor an\rfloor\leq k\leq l\leq n}\frac{|\tilde X_k-\tilde X_l|}{\sqrt{\chi_n(\tau)}}\geq \epsilon\right)=0.
\end{equation}
Since $\sup_{\lfloor an\rfloor\leq k\leq l\leq n}|\tilde X_k-\tilde X_1|\leq 2\sup_{\lfloor an\rfloor\leq s\leq n}|\tilde X_s-\tilde X_{\lfloor an\rfloor}|$, we can deal with $|\tilde X_s-\tilde X_{\lfloor an\rfloor}|$ in the following proof. Let $m$ be $|\tilde X_{m}-\tilde X_{\lfloor {a}n\rfloor}|=\sup_{\lfloor an\rfloor\leq s\leq n}|\tilde X_s-\tilde X_{\lfloor an\rfloor}|$. the decomposition of it is
\begin{equation}\label{decomposition2}
\tilde X_{m}-\tilde X_{\lfloor an\rfloor}=\sum_{j=\tau(k(\lfloor an\rfloor)+1)}^{\tau(k(m))}\tilde Y_j+\tilde{\bar{Y}}_{m}-\tilde{\bar{Y}}_{\lfloor an\rfloor},
\end{equation}
or just $\tilde{\bar{Y}}_{m}-\tilde{\bar{Y}}_{\lfloor an\rfloor}$ if $k(\lfloor an\rfloor)=k(m)$. Follow what we did in the proof of slow cooling,
\begin{equation}
      \mathbb{P}\left(\frac{|\sum_{j=\tau(k(\lfloor an\rfloor)+1)}^{\tau(k(m))}\tilde Y_j|}{\sqrt{\chi_n(\tau)}}\geq \frac{\epsilon}{2}\right)\leq
      \mathbb{P}\left(\sup_{l\in [0,k(n)-k(\lfloor an\rfloor)]}\frac{|M_l|}{\sqrt{\chi_n(\tau)}}\geq \frac{\epsilon}{2}\right).
\end{equation}
Combining \eqref{doob} and \eqref{bound} under the case $q=1$, recall that $k(n)-k(\lfloor an\rfloor)\sim -(1/C)\log a$, there exists $C_2>0$, depending only on $\epsilon$,
\begin{equation}\label{27}
      \mathbb{P}\left(\frac{|\sum_{j=\tau(k(\lfloor an\rfloor)+1)}^{\tau(k(m))}\tilde Y_j|}{\sqrt{\chi_n(\tau)}}\geq \frac{\epsilon}{2}\right)\leq \frac{C_2}{\log n}.
\end{equation}
For $\tilde{\bar{Y}}_{m}$, following all the steps from \eqref{17} to \eqref{21}, there exists $C''>0$, depending only on $\epsilon$,
\begin{equation}\label{28}
      \mathbb{P}\left(\frac{|\tilde{\bar{Y}}_{m}|}{\sqrt{\chi_n(\tau)}}\geq \frac{\epsilon}{4}\right)\leq \frac{C''}{\log^2n}\left({k(n)-k(\lfloor an\rfloor)}\right),
\end{equation}
and obviously the right hand side goes to zero as n goes to infinity. By \eqref{26}, \eqref{27} and \eqref{28},
\begin{equation}
\begin{split}
&\mathbb{P}\left(\sup_{\lfloor an\rfloor\leq s\leq n}\frac{|\tilde X_s-\tilde X_{\lfloor an\rfloor}|}{\sqrt{\chi_n(\tau)}}\geq \epsilon\right)=
\mathbb{P}\left(\frac{|\tilde X_m-\tilde X_{\lfloor an\rfloor}|}{\sqrt{\chi_n(\tau)}}\geq \epsilon\right)\\
\leq & \mathbb{P}\left(\frac{|\sum_{j=\tau(k(\lfloor an\rfloor)+1)}^{\tau(k(m))}\tilde Y_j|}{\sqrt{\chi_n(\tau)}}\geq \frac{\epsilon}{2}\right)+
\mathbb{P}\left(\frac{|\tilde{\bar{Y}}_{m}|}{\sqrt{\chi_n(\tau)}}\geq \frac{\epsilon}{4}\right)+
\mathbb{P}\left(\frac{|\tilde{\bar{Y}}_{\lfloor an\rfloor}|}{\sqrt{\chi_n(\tau)}}\geq \frac{\epsilon}{4}\right)\\
\leq & \frac{C_2}{\log n}+\frac{C''}{\log^2n}\left({k(n)-k(\lfloor an\rfloor)}\right)+\frac{16\mathbb{V}ar (\bar Y_{\lfloor an\rfloor})}{\epsilon^2\chi_n(\tau)}\\
= & O(\frac{1}{\log n}).
\end{split}
\end{equation}
The tightness condition holds. Hence \eqref{functional limit exponential} is proved.   \qed

\bibliography{biblio}

\end{document}